\documentclass[10pt]{article}
\usepackage{graphics}
\usepackage{amssymb}
\usepackage{amsthm}

\newtheorem{theorem}{Theorem}[section]

\newtheorem{lemma}{Lemma}[section]
\newtheorem{corollary}{Corollary}[section]

\def\qed{\hfill \rule{2.5mm}{2.5mm}}
\date{}
\begin{document}
\title{On existence of $[a,b]$-factors avoiding given subgraphs
\thanks{This work is supported by Shandong Provincial Taishan Scholarship of China
and Natural Sciences and Engineering Research Council of Canada.}}
\author{ Yinghong Ma$^1$$^2$ and Qinglin Yu$^1$$^3$
\\ {\small  $^1$Center for Combinatorics, LPMC, Nankai University}
\\ {\small Tianjing, China}
\\ {\small $^2$ School of Management }
\\ {\small Shandong Normal University, Jinan, Shandong, China}
\\ {\small  $^3$Department of Mathematics and Statistics}
\\ {\small Thompson Rivers University, Kamloops, BC, Canada}
}

\maketitle

\begin{center}
\begin{minipage}{114mm}
\begin{center}{\bf Abstract}\end{center}
For a graph $G = (V(G), E(G))$, let $i(G)$ be the number of isolated
vertices in $G$. The {\it isolated toughness} of $G$ is defined as
$I(G) = min\{|S|/i(G-S) :  S\subseteq V(G), i(G-S)\geq 2\}$ if $G$
is not complete; $I(G)=|V(G)|-1$ otherwise.  In this paper, several
sufficient conditions in terms of isolated toughness are obtained
for the existence of $[a, b]$-factors avoiding given subgraphs,
e.g., a set of vertices, a set of edges and a matching,
respectively.

\end{minipage}
\end{center}


{\bf AMS(1991) Subject Classification:} 05C70

\vskip 0.5cm

\section{Introduction}

Let $G = (V(G), E(G))$ be a graph with vertex set $V(G)$ and edge set $E(G)$.
 We use $d_G(x)$ to denote the degree of $x$ in $G$ and
$\delta(G)$ to denote the minimum vertex degree of $G$. For a vertex
set $S \subseteq V(G)$, the subgraph of $G$ induced by $S$ is
denoted by $G[S]$, $i(G-S)$ and $c(G-S)$ are used for the number of
isolated vertices and the number of components in $G-S$,
respectively. A subset $I$ of $V(G)$ is  {\it an independent set} if
no two vertices of $I$ are adjacent in $G$ and a set $C$ of $V(G)$
is {\it a covering set} if every edge of $G$ is incident to a vertex
in $C$. For any two subsets $S,T\subseteq V(G)$, $E(S,T)=\{uv\in
E(G):u\in S, v\in T\}$.

Let $H$ be a spanning subgraph of $G$ and $a$, $b$ be two
nonnegative integers satisfying $a \leq b$. We call $H$ an
$[a,b]$-{\it factor} of $G$ if $a\leq d_H(x)\leq b$ for each $x\in
V(G)$. When $a=1$ and $b=m>1$, it is not hard to see that existence
of $[1,m]$-factor is equivalent to the existence of a spanning
subgraph consisting of stars with no more than $m$ edges.  So $[1,
m]$-factors are also referred as {\it star-factors}, denoted by
$S(m)$-factor.  For $a=b=k >0$, $[a, b]$-factor is commonly known as
$k$-{\it factor}. In particular, $1$-factors are often referred as
perfect matchings.

  Matching problem as one of most well-established branches of graph
theory, does not only lie at the heart of many applications, it also
gives rise to some most matured techniques (e.g., augmenting path)
and generates some deep mathematical discoveries (e.g., matching
polytope theory). Since the characterization of perfect matchings
were given by Tutte in 1947, the concept of perfect matching has
been extended to several general forms, from $k$-factors to
$f$-factors, to $[a, b]$-factors, to $(g, f)$-factors. In this
paper, we use a new graphic parameter -- isolated toughness -- to
establish several sufficient conditions for the existence of $[a,
b]$-factors with given properties. In particular, we studied the
existence of $[a, b]$-factors avoiding a set of vertices, a set of
edges and a matching, respectively.

    The new parameter, isolated toughness,  is motivated by
Chv\'{a}tal's celebrated graphic parameter, toughness.  It can be
obtained from the definition of toughness by replacing $c(G-S)$ by
$i(G-S)$. The {\it isolated toughness} $I(G)$ was first introduced
by Ma and Liu \cite{Ma} and is defined as

$$I(G)=\left\{\begin{array}{ll} min \{\frac{|S|}{i(G-S)} \ : \ S \subseteq
V(G),i(G-S)\geq 2\}& \mbox{if $G$ is not complete}; \\ |V(G)|-1 &
\mbox{otherwise}.\end{array}\right.
$$

    To study the existence of $[a, b]$-factors, we will use a necessary
and sufficient condition of $(g<f)$-factors given by  Heinrich {\it
et al.} \cite{Heinrich}.

\begin{theorem} (Heinrich et al. \cite{Heinrich})\label{theorem1}
Let $g(x)$ and $f(x)$ be nonnegative integral-valued
 functions defined on $V(G)$. If either one of the following conditions holds

    (i) $g(x)<f(x)$ for every vertex $x \in V(G)$;

    (ii) $G$ is bipartite;

\noindent then $G$ has a $(g,f)$-factor if and only if for any set
$S$ of $V(G)$ $$g(T)-d_{G-S}(T)\leq f(S)$$
 where $T=\{x : x\in V(G)-S,\,\,d_{G-S}(x)\leq g(x) \}$.
\end{theorem}

    In the above theorem, to confirm a graph possessing $(g, f)$-factors,
we need only to verify the much simpler inequality above for every
vertex set $S$, in contrast with the verification of a more complex
inequality for all possible pair of disjoint vertex sets $(S, T)$ in
Lov\'{a}sz's original characterization of general $(g, f)$-factors.
This simpler criterion enables us to deal with factor problems with
additional properties.

Let $g(x)=a<b=f(x)$ in Theorem \ref{theorem1}, it yields
a necessary and sufficient condition for existence of
$[a,b]$-factors. If $a=1$ and $b=m\geq 2$, then it becomes the
necessary and sufficient condition for a graph having
$S(m)$-factors.

\begin{theorem} (Anstee \cite{Anstee})\label{theorem2}
Let $G$ be a graph and let $a<b$ be two positive integers. Then
$G$ has an $[a,b]$-factor if and only if for any
 $S\subseteq V(G)$, $$a|T|-d_{G-S}(T)\leq b|S|$$ holds,
 where $T=\{x:x\in V(G)-S,\,\, d_{G-S}(x)\leq a-1\}$.
\end{theorem}

\noindent {\bf Remarks:} Let $T=\{x:x\in V(G)-S,\,\, d_{G-S}(x)\leq
a-1\}$, $T'=\{x : x\in V(G)-S,\,\, d_{G-S}(x)\leq a\}$ and $T''=\{x
: x\in V(G)-S,\,\,d_{G-S}(x)=a\}$. Then $T = T' - T''$. Since
$a|T''|-d_{G-S}(T'')=0$, we have
$a|T|-d_{G-S}(T)=(a|T'|-d_{G-S}(T')) - (a|T''|-d_{G-S}(T''))$ =
$a|T'|-d_{G-S}(T')$.
So $T$ in Theorem \ref{theorem2} is equivalent to that in
Theorem \ref{theorem1} when $g(x)=a$.\\

    Use the isolated toughness as a sufficient condition, Ma and
Liu \cite{Ma} provided an existence theorem for $[a, b]$-factors.

\begin{theorem} (Ma and Liu \cite{Ma})\label{theorem3}
Let $G$ be a graph with $\delta (G)\geq a$ and $I(G)\geq
a-1+\frac{a}{b}$. Then $G$ has $[a,b]$-factors.
\end{theorem}

    For convenience, we denote $\delta_G(a,b;S)=b|S|-a|T|+d_{G-S}(T)$.
So Theorem \ref{theorem2} can be restated as that $G$ has
$[a,b]$-factors if and only if $\delta_G(a,b;S)\geq 0$ for any
$S\subseteq V(G)$.

\section{Main Results}
Throughout the paper, we always assume that $a,b,m$ and $n$ are
positive integers satisfying $1 \leq a<b$. So we will not reiterate
these conditions again in the theorems or proofs.

    The first result is to investigate the existence of $[a,b]$-factors in the operation of
vertex-deletion.

\begin{theorem} \label{theoremA}
Let $G$ be a graph with $\delta (G)\geq a+n$ and
the isolated toughness $I(G)\geq a-1+n+\frac{a-1}{b}$. Then, for any
$n$-subset $V'\subset V(G)$, $G-V'$ has $[a,b]$-factors.
\end{theorem}

The condition $I(G)\geq a-1+n+\frac{a-1}{b}$ in Theorem
\ref{theoremA} can not be weakened, that is, if we replace the
condition by $I(G)\geq a-1+n+\frac{a-1}{b}-\epsilon$, where
$\epsilon$ is any positive real number, then there exists an $n$-set
$V_0\subset V(G)$ such that $G-V_0$ has no $[a,b]$-factor. Consider the
following family of graphs.

Construct $H$ as follows:

$V(H)=V(K_{m(a-1)})\cup V((mb+1)K_1)\cup V(K_{(mb+1)(a-1+n)}),$

$E(H)=E(K_{m(a-1)})\cup E(K_{(mb+1)(a-1+n)}) \cup
(\cup_{i=1}^{mb+1}u_iv_i) \cup \{xy : x\in V(K_{m(a-1)})$, $y\in
(mb+1)K_1\}$,

\noindent where $V((mb+1)K_1)=\{v_1,v_2,\cdots,v_{mb+1}\}$ and
$\{u_1,u_2,\cdots,u_{mb+1}\}\subset V(K_{(mb+1)(a-1+n)})$. Let
$S=V(K_{m(a-1)})\cup V(K_{(mb+1)(a-1+n)}).$ Clearly, $I(H)\leq
\frac{|S|}{i(G-S)} = \frac{(mb+1)(a-1+n)+m(a-1)}{mb+1}\rightarrow
a-1+n+\frac{a-1}{b}$ when $m\rightarrow +\infty$, and is less than
$a-1+n+\frac{a-1}{b}$. Let $V_0\subset V(K_{(mb+1)(a-1+n)})
\backslash \{u_1,\cdots,u_{mb+1}\}$ be an $n$-vertex set, then
$H-V_0$ has no $[a,b]$-factors.  To see this, consider the set
$S=V(K_{m(a-1)})\subset V(H)-V_0$, then we have $T=V((mb+1)K_1)$ and
$a|T|-d_{H-V_0-S}(T)=(mb+1)(a-1)>mb(a-1)=b|S|$. Thus, by Theorem
\ref{theorem2}, $H-V_0$ has no $[a,b]$-factor. So in this sense
 Theorem \ref{theoremA} is best possible.

For the existence of $[a,b]$-factors resulting from the operation of
edge-deletion, we first investigate star-factors and obtain the
following.

\begin{theorem} \label{theoremB}
Let $G$ be a graph with $\delta (G)\geq 1+n$ and $I(G)\geq
\frac{1}{m-n}$,  where  $1\leq n\leq \frac{m}{2}$. Then for any $n$-subset $E'\subset E(G)$,
$G-E'$ has $S(m)$-factors.
\end{theorem}

    A sufficient condition for the existence of $[a, b]$-factors in the operation of
matching-deletion is given below.

\begin{theorem}\label{theoremC}
If a graph $G$ satisfies $\delta (G)\geq a+n$ and $I(G)\geq
a-1+\frac{a+2n-1}{b}$, then for any $n$-matching $M$ of $G$, $G-M$
has $[a,b]$-factors.
\end{theorem}

We next investigate hierarchy relation for the operation of vertex-deletion.

\begin{theorem} \label{theoremD}
Let $G$ be a graph with $\delta (G)\geq a+n$. If, for any arbitrary
$n$-subset $V'\subset V(G)$, $G-V'$ has $[a,b]$-factors, then, for
any $(n-1)$-subset $V''\subset V(G)$,  $G-V''$ has $[a,b]$-factors
as well.
\end{theorem}

Finally we present a different type of sufficient condition for the
existence of $[a,b]$-factors excluding any edge of  $E(G)$.

\begin{theorem} \label{theoremE}
Let $G$ be a graph with $\delta (G)\geq a+2$. If $G-\{x,y\}$ has
$[a,b]$-factors for every pair of vertices $x,y\in V(G)$, then $G-e$
has $[a,b]$-factors for any given edge $e\in E(G)$.
\end{theorem}

\section{Proofs of Theorems \ref{theoremA} and \ref{theoremC}}

    In order to prove Theorem \ref{theoremA}, we need the following
lemmas.

\begin{lemma} \label{lemmac}
Let $G$ be a graph. Then, for any $n$-subset $V'\subset V(G)$,
$G-V'$ has an $[a,b]$-factor if and only if for any $S\subset V(G)$
with $V'\subseteq S$
$$\delta_{G}(a,b;S)=b|S|-a|T|+d_{G-S}(T)\geq bn$$
where $T=\{x : x\in V(G)-S, \, d_{G-S}(x)\leq a-1\}$.
\end{lemma}
\noindent {\it\bf Proof: }  Suppose that for any $n$-subset
$V'\subset V(G)$, $G-V'$ has $[a,b]$-factors. Let $G'=G-V'$, by
Theorem \ref{theorem2}, $G'$ has $[a,b]$-factor if and only if for
any $S'\subset V(G')$, $\delta_{G'}(a,b;S')\geq 0$, where $T'=\{x :
x\in G'-S', d_{G'-S'}(x)\leq a-1\}$. Let $S'=S - V'$, then $T=T'$
and
$\delta_{G'}(a,b;S')=b|S-V'|-a|T|+d_{G'-S'}(T)=\delta_{G}(a,b;S)-bn$.
Therefore, $\delta_{G}(a,b;S)\geq bn$ since $G'-S'=G-S$ and
$\delta_{G'}(a,b;S')\geq 0$.

Conversely, suppose there exists some $n$-subset $V_0\subset V(G)$
such that $G'=G-V_0$ has no $[a,b]$-factor. By Theorem
\ref{theorem2}, there exists $S_0\subset V(G')$ such that
$\delta_{G'}(a,b;S_0)<0$, where $T_0=\{x : x\in G'-S_0,
d_{G'-S_0}(x)\leq a-1\}$. Let $S=S_0\cup V_0$. Then $G'-S_0=G-S$ and
$T=T_0$, and thus

$$\begin{array}{rcl}
\baselineskip 25pt \delta_{G}(a,b;S)& = & b|S_0\cup
V_0|-a|T_0|+d_{G'-S_0}(T_0)\\
&=&b|S_0|-a|T_0|+d_{G'-S_0}(T_0)+b|V_0|\\
&=&\delta_{G'}(a,b;S_0)+bn,
\end{array}$$

\noindent therefore, $\delta_{G}(a,b;S) < bn$, a contradiction.
Hence, $G-V_0$
has $[a,b]$-factors for any $n$-subset $V_0\subset V(G)$. \qed \\

    To prove the main lemma (Lemma \ref{lemmaD1}), we will require a
technical tool here stated as a corollary below which is an enriched
version of the following result from Katerinis \cite{Katerinis}.

\begin{lemma} (Katerinis \cite{Katerinis})\label{lemmaD}
Let $H$ be a graph and $S_1,S_2,\cdots, S_{a-1}$ a vertex partition
of $H$ such that $d_H(x)\leq j$ for each $x\in S_j\,(1\leq j\leq
a-1)$. Then there exist an independent set $I$ and a covering set
$C$ of $H$ such that
$$\sum_{j=1}^{a-1}(a-j)c_j\leq \sum_{j=1}^{a-1}j(a-j)i_j,$$ where
$c_j=|S_j\cap C|$ and $i_j=|S_j\cap I|$.
\end{lemma}

\begin{corollary} \label{coroA}
Let $H$ be a graph and $S_1,S_2,\cdots, S_{a-1}$ a vertex partition
of $H$ such that $d_H(x)\leq j$ for each $x\in S_j\,(1\leq j\leq
a-1)$. Then there exist a {\bf maximal} independent set $I$ and a
covering set $C$ of $H$ such that $I \cap C = \emptyset$ and
$$\sum_{j=1}^{a-1}(a-j)c_j\leq \sum_{j=1}^{a-1}j(a-j)i_j,$$ where
$c_j=|S_j\cap C|$ and $i_j=|S_j\cap I|$.
\end{corollary}
\noindent {\it\bf Proof: }  From Lemma \ref{lemmaD}, there exist an
independent set $I'$ and a covering set $C'$ of $H$ such that
$$\sum_{j=1}^{a-1}(a-j)c_j'\leq \sum_{j=1}^{a-1}j(a-j)i_j',$$ where
$c_j'=|S_j\cap C'|$ and $i_j'=|S_j\cap I'|$.

    Note the fact that any complement of an independent set must be
a covering set.  Let $I$ be a {\it maximal} independent set
containing $I'$, $C = V(G) - I$ and $C'' = C' - (C' \cap I')$. Then
$C$ and $C''$ are both covering sets. Thus $I \supseteq I'$, $C
\subseteq C'' \subseteq C'$ and $I \cap C = \emptyset$. Since
$c_j=|S_j\cap C| \leq |S_j\cap C'| =c_j'$ and $i_j'=|S_j\cap I'|
\leq |S_j\cap I| = i_j$, we have
$$\sum_{j=1}^{a-1}(a-j)c_j\leq \sum_{j=1}^{a-1}(a-j)c_j'\leq
\sum_{j=1}^{a-1}j(a-j)i_j' \leq \sum_{j=1}^{a-1}j(a-j)i_j.$$
\qed \\

    The techniques used to prove Theorems \ref{theoremA} and \ref{theoremC} are
along the same line, so we present the main ideas as a lemma below.

\begin{lemma} \label{lemmaD1}
Let $2 \leq k \leq b$. If a graph $G$ satisfies $\delta (G)\geq a+n$
and $I(G)\geq a-1+\frac{a+kn-1}{b}$, then
$\delta_{G}(a,b;S)=b|S|-a|T|+d_{G-S}(T)\geq kn$ for any subset $S
\subseteq V(G)$ with $T =\{x : x\in V(G)-S, \, d_{G-S}(x)\leq
a-1\}\not= \emptyset$.
\end{lemma}
\noindent {\it\bf Proof: }
Use the argument of contradiction.  Suppose that there exists a vertex set
$S \subseteq V(G)$ such that
$$\delta_G(a,b;S)=b|S|-a|T|+d_{G-S}(T)<kn,\eqno(1)$$
where $T=\{x : x\in V(G)-S, d_{G-S}(x)\leq a-1\}\not= \emptyset$.

For each $0\leq j\leq a-1$, let $T^j=\{x : x\in T, d_{G-S}(x)=j\}$
($T^j$ may be an empty set) and $|T^j|=t_j$. Let $H=G[T^1\cup
T^2\cup \cdots \cup T^{a-1}]$, clearly $\{T^1,T^2,\cdots,T^{a-1}\}$
is a vertex partition of $H$ and $d_H(x)\leq j$ for each $x\in
T^j$\,($1\leq j\leq a-1$). Then, by Corollary \ref{coroA}, there
exist a {\it maximal} independent set $I$ and a covering set $C$ of
$H$ such that $I \cap C= \emptyset$ and
$$\sum_{j=1}^{a-1}(a-j)c_j\leq \sum_{j=1}^{a-1}j(a-j)i_j,\eqno (2)$$
where $c_j=|T^j\cap C|$ and $i_j=|T^j\cap I|$, $j=1,2,\cdots,a-1$.

Let $W=G-(S\cup T)$ and $U=S\cup C\cup (N_{G-S}(I)\cap V(W))$, we have
$$|U|\leq |S|+\sum_{j=1}^{a-1}ji_j,\eqno (3)$$
and
$$i(G-U) \geq t_0 +|I| = t_0+\sum_{j=1}^{a-1}i_j.\eqno(4)$$

{\it Case 1.} $t_0+\sum_{j=1}^{a-1}i_j\leq 1$.

    Since $T \not= \emptyset$, it follows either $t_0=1$ and
$\sum_{j=1}^{a-1}i_j=0$ or $t_0=0$ and $\sum_{j=1}^{a-1}i_j=1$.

    If $t_0=1$ and $\sum_{j=1}^{a-1}i_j=0$, then $H=\emptyset$.
Let $T=\{v\}$, by (1), we have $a+kn>b|S|\geq
b(a+n) \geq a+kn$ as $|S|\geq d_G(v)\geq \delta(G)\geq a+n$ and $b \geq k$,
a contradiction.

    If $t_0=0$ and $\sum_{j=1}^{a-1}i_j=1$, then, for
some $j_0\in \{1,2,\cdots,a-1\}$, $i_{j_0}=1$ and $i_j=0$ for all
$j\in \{1,2,\cdots,a-1\}-{j_0}$.  Let
$I=\{u\}$, then $a+n\leq \delta(G)\leq d_G(u)\leq |S|+j_0$ or
$|S|\geq a+n-j_0$. Therefore,
$$b|S|-kn\geq b(a+n-j_0)-kn = b(a-j_0) +(bn-kn). \eqno(5)$$

Since $I$ is maximal, we see $V(H) \subseteq I \cup C$ and thus $t_j \leq i_j + c_j$.
Recall that $t_0 = 0$, by (2), it yields $a|T|-d_{G-S}(T)=\sum_{j=1}^{a-1}(a-j)t_j +at_0 \leq
\sum_{j=1}^{a-1}(a-j)i_j+\sum_{j=1}^{a-1}(a-j)c_j
 \leq a-j_0+\sum_{j=1}^{a-1}j(a-j)i_j = a-j_0+j_0(a-j_0)$.
Combining (1),  (5) and the previous inequality,  we have
$$b(a-j_0)+(bn-kn) \leq b|S| - kn < a-j_0 + j_0(a-j_0)$$

\noindent or
$$ba-a<-j_0^2+aj_0+bj_0-j_0. \eqno(6)$$

Let $f(x)=-x^2+(a+b-1)x$. Then the maximum value of the quadratic function $f(x)$ is
$\frac{(a+b-1)^2}{4}$ when $x = \frac{a+b-1}{2}$. However, $f(x)$ can not attain this value
since $x\in \{1,2,\cdots,a-1\}$.  Because $f(1)<f(2)<\cdots<f(a-1)=b(a-1)$, (6) becomes
$ba-a<-j_0^2+aj_0+bj_0-j_0 \leq f(a-1) = ba-b$,  a contradiction.\\

    {\it Case 2}. $t_0+\sum_{j=1}^{a-1}i_j\geq 2$.

From (4), we have $i(G-U)\geq t_0+\sum_{j=1}^{a-1}i_j \geq 2$. By the definition of $I(G)$
and (4),  we have
$$|U|\geq I(G)i(G-U)\geq (t_0+\sum_{j=1}^{a-1}i_j)I(G),$$
\noindent or
$$|S|\geq \sum_{j=1}^{a-1}(I(G)-j)i_j+t_0I(G).\eqno (7)$$
Recall $t_j \leq i_j + c_j$, thus (1), (2) and (7) imply

$$\begin{array}{lcll}
    a|T|-d_{G-S}(T)&= &\sum_{j=1}^{a-1}(a-j)t_j+at_0\\
                   &\leq &\sum_{j=1}^{a-1}(a-j)i_j+\sum_{j=1}^{a-1}(a-j)c_j+at_0 \\
    & \leq & \sum_{j=1}^{a-1}(a-j)i_j+\sum_{j=1}^{a-1}j(a-j)i_j+at_0,
\end{array}$$

and
$$a|T|-d_{G-S}(T)
> b|S| -kn \geq \sum_{j=1}^{a-1}(bI(G)-bj)i_j+bt_0I(G)-kn.\eqno(8)$$
Therefore,
$$\sum_{j=1}^{a-1}(-j^2+(a+b-1)j)i_j>\sum_{j=1}^{a-1}(bI(G)-a)i_j+bt_0I(G)-at_0-kn.\eqno (9)$$

If $\sum_{j=1}^{a-1}i_j=0$, then $t_0\geq 2$, $H=\emptyset$ and
$|T|=t_0=i(G-S)$. So $\frac{|S|}{i(G-S)}\geq I(G)$ implies $|S|\geq
I(G)t_0$. By (1), $b|S|<at_0+kn$ and thus $at_0+kn>b|S|\geq
bt_0I(G)\geq at_0+kn$ since $I(G)\geq a-1+\frac{a+kn-1}{b}$, a
contradiction.

 If $\sum_{j=1}^{a-1}i_j\not=0$,
 then we can see $bt_0I(G)-at_0-kn = t_0(bI(G)-a)-kn \geq (1-kn)\sum_{j=1}^{a-1}i_j$ by
noting $bI(G)-a \geq 0$ and recalling that $t_0+\sum_{j=1}^{a-1}i_j\geq 2$ and $\sum_{j=1}^{a-1}i_j\not=0$.
 From (9), we obtain $$\sum_{j=1}^{a-1}(-j^2+(a+b-1)j)i_j>\sum_{j=1}^{a-1}(bI(G)-a+1-kn)i_j.$$
Therefore, there is at least one $j\in \{1,2,\cdots,a-1\}$ such that
$-j^2+(a+b-1)j>bI(G)-a+1-kn$. But this is impossible, because
$-j^2+(a+b-1)j\leq b(a-1)$ for all the $j\in \{1,2,\cdots,a-1\}$ and
$bI(G)-a+1-kn\geq b(a-1)$ as $I(G) \geq a-1+\frac{a+kn-1}{b}$.

  The lemma is proven. \qed

\vspace{3mm}
    With Lemma \ref{lemmaD1} in the hand, we can provide short proofs for Theorems \ref{theoremA} and
\ref{theoremC}. \\

\noindent {\it\bf Proof of Theorem \ref{theoremA}:} If $G$ is a
complete graph, clearly the theorem holds. So we assume that $G$ is
not complete.

Suppose that $G$ satisfies the conditions of the theorem, but there
exists an $n$-subset $V_0\subset V(G)$ such that $G'=G-V_0$ has no
$[a,b]$-factor. By Lemma \ref{lemmac}, there exists a vertex set $S$
with $V_0\subset S$ such that
$$\delta_G(a,b;S)=b|S|-a|T|+d_{G-S}(T)<bn,\eqno(10)$$
where $T=\{x : x\in V(G)-S, d_{G-S}(x)\leq a-1\}$.

    If $T=\emptyset$, then (10) becomes $bn > \delta_G(a,b;S)=b|S|
\geq bn$ as $|S| \geq n$, a contradiction.

     If $T \not= \emptyset$, applying Lemma \ref{lemmaD1} with $k=b$ we conclude
that (10) does not hold.

    So we conclude that $G-V_0$ has $[a,b]$-factors for
any $n$-subset $V_0\subset V(G)$.  \qed \\

Next, we consider the existence of $[a,b]$-factors excluding an
$n$-matching. \\

\noindent{\bf Proof Theorem \ref{theoremC}:} Suppose that $G$
satisfies the conditions given in the theorem, but there
exists a matching $M$ in $G$ with $|M|=n$ such that $G-M=G'$ has no
$[a,b]$-factor. By Theorem \ref{theorem2}, there exists some $S\subset
V(G')=V(G)$ such that
$$a|T'|-d_{G'-S}(T')>b|S|\eqno(11)$$ where
$T'=\{x : x\in V(G')-S, d_{G'-S}(x)\leq a-1\}$. Denote $T=\{x : x\in
V(G)-S, d_{G-S}(x)\leq a-1\}$.

Clearly, $S\not=\emptyset$. Otherwise, $T'=\emptyset$ since
$\delta(G)\geq a+n$ and then, by (11), $a|T'|-d_{G'-S}(T')=0>b|S|=0$,
a contradiction.

If $V(M)\subseteq S$, then $T = T'$ and $d_{G'-S}(T')=d_{G-S}(T)$.
Since $\delta (G)\geq a+n$ and $I(G)\geq a-1+\frac{a+2n-1}{b}\geq
(a-1)+\frac{a}{b}$, by Theorem \ref{theorem3}, $G$ has
$[a,b]$-factors or $b|S|\geq
a|T|-d_{G-S}(T)=a|T'|-d_{G'-S}(T')>b|S|$, a contradiction to (11).
So we assume $V(M)\not\subseteq S$.

    Let $W=G-(S\cup T)$ and $V_0=V(M)$. Denote $V_W= \{ x\in V_0\cap W:
d_{G-S}(x)=a$ and $ \exists y \in G-S$ so that $xy\in M\}$.
Clearly, $T'=T\cup V_W$ and the degrees of vertices of $V_W$ in $G'-S$ are $a-1$.
Therefore, $d_{G'-S}(T')=d_{G'-S}(T)+d_{G'-S}(V_W)\geq
d_{G-S}(T)+d_{G-S}(V_W)-2n$ and $d_{G-S}(V_W)=a|V_W|$. By (11),
$b|S|<a|T'|-d_{G'-S}(T')=a|T|+a|V_W|-d_{G'-S}(T')\leq
a|T|-d_{G-S}(T)+2n$.

    From the above discussion, to prove the theorem we need only to show that
the following inequality does not hold for any $S\subset V(G)$
$$b|S|-a|T|+d_{G-S}(T)<2n.\eqno(12)$$

    For any $S\subset V(G)$, if $T=\emptyset$, from (11), $T' \not= \emptyset$ and
thus there exists a vertex $u \in T'$
so that $d_{G-S}(u) = a$. Thus $|S| \geq n$ as $\delta(G) \geq a+n$.
So (12) becomes $0 > b|S| - 2n \geq bn-2n \geq 0$, that is, (12) does not hold.

     If $T \not= \emptyset$, applying Lemma \ref{lemmaD1} with $k=2$ we conclude
that (12) does not hold.

 We complete the proof. \qed

\section{Proofs of Theorems \ref{theoremB}, \ref{theoremD} and \ref{theoremE}}

In order to prove Theorem \ref{theoremB}, we need the following
lemmas.

\begin{lemma} (Las Vergnas \cite{Las})\label{lemmaE}
Let $G$ be a graph. Then $G$ has $S(m)$-factors if and only if
 $i(G-S)\leq m|S|$ for any $S\subset V(G)$.
\end{lemma}

    Lemma \ref{lemmaE} can be derived from Theorem \ref{theorem1} easily by letting
$a=1$ and $b=m>1$ . Using the notation of isolated toughness,
Lemma \ref{lemmaE} can be restated as that $G$ has $S(m)$-factors if
and only if $I(G)\geq \frac{1}{m}$.

The following proposition can be seen easily, so we omit the proof.

\begin{lemma}\label{lemmaF}
For any edge $e$ of a graph $G$, then $i(G)\leq i(G-e)\leq
i(G)+2$.
\end{lemma}

\vspace{3mm}
 Now we turn to the proof of Theorem \ref{theoremB}.\\

\noindent {\bf Proof of Theorem \ref{theoremB}:}
 Let $G$ be a graph satisfying the conditions given in the theorem, but there
exists an edge set $E_0\subset E(G)$ with $|E_0|=n\leq \frac{m}{2}$
such that $G-E_0$ has no $S(m)$-factor. Setting $G-E_0=G'$, then, by
Lemma \ref{lemmaE}, $I(G')<\frac{1}{m}$. That is, there exists a
vertex set $S\subset V(G')=V(G)$ such that

$$i(G'-S)>m|S|.\eqno(12)$$
Clearly, $S\not=\emptyset$ (since $\delta(G)\geq 1+n$).

 By Lemma \ref{lemmaF}, $i(G'-S) = i(G-E_0-S) \leq i(G-S)+2n$.
We consider the following cases.\\

{\it Case 1.} $i(G-S)\geq 2$. Then, by the definition of $I(G)$, we
have $i(G'-S)\leq i(G-S)+2n\leq (m-n)|S|+2n$ since $I(G)\geq
\frac{1}{m-n}$.

If $|S|\geq 2$, then $i(G'-S)\leq i(G-S)+2n\leq (m-n)|S|+2n\leq
m|S|$, a contradiction to (12).

If $|S|=1$, let $u,v$ be two isolated vertices in $G-S$, then
$d_{G}(u)=d_G(v)=1$ since $S$ is a cut set of $u$ and $v$, but this
is impossible since $\delta(G)\geq 1+n>1$.

{\it Case 2.} $i(G-S)=0$. In this case, $m\leq m|S|<i(G'-S)\leq 2n,$
a contradiction to the condition $n\leq \frac{m}{2}$.

{\it Case 3.} $i(G-S)=1$. Then $|S|\geq n+1$ and thus $2n+2\leq m(n+1)\leq m|S|<i(G'-S)\leq
i(G-S)+2n=2n+1,$ a contradiction.

Therefore, $G-E_0$ has $S(m)$-factors for any $n$-subset
$E_0\subset E(G)$.\qed \\

\noindent{\it\bf Proof Theorem \ref{theoremD}:}
We verify the theorem for the case of $n = 1$ first, i.e., the following claim:

{\it Claim.} If $G-x$ has $[a,b]$-factors for any $x\in V(G)$, then
$G$ has $[a,b]$-factors.

    Otherwise, $G$ has no $[a,b]$-factors and thus, by Theorem \ref{theorem2},
there exists $U\subset V(G)$ such that $a|W|-d_{G-U}(W)>b|U|$, where
$W=\{x: x\in V(G)-U, d_{G-U}(x)\leq a-1\}$.  Choose a vertex $v$ from
$U$, let $U'=U-\{v\}$, then $(G-v)-U'=G-U$ and $\{ x: x\in
V(G-v)-U', \ \ d_{(G-v)-U'}(x)\leq a-1\}=W$. Therefore we have
$a|W|-d_{(G-v)-U'}(W)\leq b|U'|=b|U|-b<b|U|$ since $G-v$ has
$[a,b]$-factors, a contradiction since $a|W|-d_{G-U}(W)>b|U|$.
Hence, $G$ has $[a,b]$-factors.

    Applying the above claim and using induction arguments, we can see that
$G-V''$ has $[a,b]$-factors for any $(n-1)$-subset $V''$
if $G-V'$ has $[a,b]$-factors for any $n$-subset $V'$. \qed \\

  Next we present a characterization for $[a, b]$-factors excluding an edge.
As an application, Theorem \ref{theoremE} can be easily derived from it.
In fact, the lemma itself is of interest.

\begin{lemma} \label{lemmaH}
Let $G$ be a graph and $e=uv$ be any edge of $G$. Then $G$ has
$[a,b]$-factors excluding the edge $e$ if and only if
$$\delta_G(a,b;S)\geq \rho(S)$$

\noindent holds for any $S\subseteq V(G)$, where $G'=G-e$, $T'=\{x :
x\in V(G)-S, d_{G'-S}(x)\leq  a-1\}$ and
$$\rho(S)= \left \{
\begin{array}{lll}
2 &  &\mbox{both $u$ and $v$ belong to $T'$}; \\
1 & &\mbox{one of $\{ u, v\}$ lies in $T'$ and the other is in $G-(S\cup T')$}; \\
0 & &\mbox{otherwise.}\end{array} \right.
$$
\end{lemma}

\noindent {\bf Proof:} Suppose that for a fixed edge $e=uv$ of
$G$, $G'=G-e$ has $[a,b]$-factors. By Theorem \ref{theorem2}, for
any $S\subset V(G')=V(G)$, $\delta_{G'}(a,b;S)\geq 0$. Let $W'=G'-(S\cup
T')$ and $T=\{x : x\in V(G)-S, d_{G-S}(x)\leq a-1\}$.

{\it Case 1.} $uv\in E(T')$. If $d_{G'-S}(u)=d_{G'-S}(v)=a-1$, then
$T=T'-\{u,v\}$,
$d_{G'-S}(T')=d_{G'-S}(T)+d_{G'-S}(\{u,v\})=d_{G-S}(T)+2(a-1)$, and
$0\leq b|S|-a|T'|+d_{G'-S}(T')=b|S|-a|T|+d_{G-S}(T)-2$ since $G'$
has $[a,b]$-factors. Therefore, $\delta_{G}(a,b;S)\geq 2.$ If
$d_{G'-S}(u)<a-1$ and $d_{G'-S}(v)<a-1$. Then $T=T'$ and
$d_{G'-S}(T')=d_{G-S}(T)-2$. Hence, $\delta_G(a,b;S)\geq 2$. If
$d_{G'-S}(u)<a-1$ and $d_{G'-S}(v)=a-1$ (or $d_{G'-S}(v)<a-1$ and
$d_{G'-S}(u)=a-1$). Then $T=T'-\{v\}$ and
$d_{G'-S}(T')=d_{G-S}(T)+a-2$. Hence, $\delta_G(a,b;S)\geq 2$.

{\it Case 2.} $uv\in E(T',W')$. Without loss of generality, let $u\in T'$ and $v\in
W'$, then we have $d_{G'-S}(u)\leq a-1$ and $d_{G'-S}(v)\geq a$. If
$d_{G'-S}(u)<a-1$, then $T=T'$. Therefore, $0\leq
\delta_{G'}(a,b;S)=b|S|-a|T'|+d_{G'-S}(T')=\delta_G(a,b;S)-1$, that
is, $\delta_G(a,b;S)\geq 1$. If $d_{G'-S}(u)=a-1$, then
$T=T'-\{u\}$. Therefore, $d_{G'-S}(T')=d_{G-S}(T)+a-1$ and then
$\delta_{G'}(a,b;S)=\delta_G(a,b;S)-1$. Hence, $\delta_G(a,b;S)\geq
1$.

{\it Case 3.} $uv\in E(S,T'\cup W')\cup E(S)\cup E(W')$. Then $T'=T$
and $d_{G'-S}(T')=d_{G-S}(T)$. Therefore, $\delta_G(a,b;S)\geq 0$.

    From the above discussion, we conclude $\delta_G(a,b;S)\geq \rho(S)$.\\

Next we prove the sufficiency. Suppose that there exists an edge
$e_0=uv\in E(G)$ such that $G'=G-e_0$ has no $[a,b]$-factor. By
Theorem \ref{theorem2}, there exists a non-empty set $S\subseteq
V(G')$ such that $\delta_{G'}(a,b;S)<0$, where $T'=\{x : x\in
V(G')-S, d_{G'-S}(x)\leq a-1\}$. Let $W'=G'-(S\cup T')$ and $T=\{x :
x\in V(G)-S, d_{G-S}(x)\leq a-1\}$.

If $e_0\in E(S,T'\cup W')\cup E(S)\cup E(W')$. Then $T=T'$ and
$d_{G'-S}(T')=d_{G-S}(T)$. Therefore,
$0>\delta_{G'}(a,b;S)=\delta_G(a,b;S)\geq 0$, a contradiction. If
$e_0\in E(T',W')$, say $u\in T'$ and $v\in W'$, we see that
$d_{G'-S}(u)\leq a-1$ and $d_{G'-S}(v)\geq a$. Then $T\subseteq T'$
and so $0>\delta_{G'}(a,b;S)=\delta_G(a,b;S)-1\geq 0$, a
contradiction. If $e_0\in E(T')$, then $T\subseteq T'$ and
$0>\delta_{G'}(a,b;S)=\delta_G(a,b;S)-2\geq 0$, a contradiction
again.

    So $G-e$ has $[a,b]$-factors for any $e\in E(G)$. \qed

\vspace{0.5cm}
\noindent {\bf Proof of Theorem \ref{theoremE}:}
Let $S$ be any subset of $V(G)$.

    If $S=\emptyset$, then $T=\emptyset$  and $\delta_G(a,b;S) = 0$.

    If $|S|=1$, then $|T|=0$ (since $\delta(G)\geq a+2$) and thus $\delta_G(a,b;S)=b|S|=b\geq 2.$

    If $|S|\geq 2$, then there exist vertices $x, y \in S$.  Let $V' = \{x, y \}$ in Lemma \ref{lemmac},
 since $G-\{x, y\}$ has $[a, b]$-factors, then we have
$\delta_G(a,b;S)\geq 2b>2$.

    Therefore, we conclude $\delta_G(a,b;S)\geq \rho(S)$ for any $S\subset V(G)$.  By Lemma \ref{lemmaH},
$G-e$ has $[a,b]$-factors. \qed

\vskip 20pt

\noindent \title{\Large\bf Acknowledgments} \maketitle
 The authors are indebted to the anonymous referees for their constructive comments.

\end {document}